\documentclass[preprint,12pt, leqno]{elsarticle}

\usepackage[margin=2.5cm]{geometry}

\usepackage{graphicx}
\usepackage{amssymb}
\usepackage{amsmath}
\usepackage{siunitx}
\usepackage{amsthm}

\usepackage{mathtools}

\usepackage{imakeidx}

\usepackage{dirtytalk}
\usepackage{extarrows}
\usepackage{accents}
\usepackage{tikz}
\usetikzlibrary{cd}
\usepackage{comment}
\usepackage{setspace}
\usepackage{enumitem}
\usepackage{epigraph}
\usepackage{soul}

\usepackage{lineno}

\usepackage[normalem]{ulem}

\journal{?}

\newcommand{\UGal}{\mathcal{G}_F}
\newcommand{\PQ}{\cP_F^+(\mQ)}

\newcommand{\mA}{\mathbb{A}}
\newcommand{\mQ}{\mathbb{Q}}

\newcommand{\mC}{\mathbb{C}}
\newcommand{\mZ}{\mathbb{Z}}

\newcommand{\cP}{\mathcal{P}}

\newcommand{\cD}{\mathcal{D}}
\newcommand{\cT}{\mathcal{T}}

\newcommand{\cC}{\mathcal{C}}
\newcommand{\cS}{\mathcal{S}}
\newcommand{\cE}{\mathcal{E}}
\newcommand{\cG}{\mathcal{G}}
\newcommand{\cU}{\mathcal{U}}

\newcommand{\fM}{\mathfrak{M}}
\newcommand{\fu}{\mathfrak{u}}

\newcommand{\SL}{\mbox{SL}}
\newcommand{\GL}{\mbox{GL}}

\newcommand{\ubf}[1]{\textit{#1}}
\newcommand{\gvec}[1]{g\mbox{Vec}_{#1}}
\newcommand{\LG}{\;^LG}
\newcommand{\whG}{\widehat{G}}
\newcommand{\Res}{\mbox{Res}}
\newcommand{\SP}{\mkern70mu }

\numberwithin{equation}{section}

\newtheorem{conjecturegroup}{Conjecture}

\makeatletter
\newcounter{subconjecture}[conjecturegroup]
\renewcommand{\thesubconjecture}{\theconjecturegroup(\alph{subconjecture})}
\newtheorem*{subconjecture*}{Conjecture} 
\newenvironment{subconjecture}[1][]
  {\refstepcounter{subconjecture}\par\vspace{0.5\baselineskip} 
   \textbf{Conjecture \thesubconjecture.} \itshape}
  {\par\vspace{0.5\baselineskip}}
\makeatother

\newtheoremstyle{named}{}{}{\itshape}{}{\bfseries}{.}{.5em}{\thmnote{#3's }#1}
\theoremstyle{named}

\begin{document}

\begin{frontmatter}


\title{Motives and Automorphic Representations}


\author{James Arthur}



\end{frontmatter}
\begin{flushright}
  To my friend and colleague, Gerard Laumon.  
\end{flushright}


\tableofcontents

\section*{Foreword}\label{sec:forword}
\addcontentsline{toc}{section}{\protect\numberline{}Foreword}

This paper is an expansion of a lecture I gave at Luminy in September 2023 with the same title. The underlying premise was the topic of an earlier note \cite{Ar2}. It was on an explicit version of the automorphic Galois group $L_F$ that would classify automorphic representations. This followed a suggestion of Langlands in \cite[Section 2]{L4} on the possibility of such a universal group that could be used to express the general analogue of the Shimura-Taniyama-Weil conjecture. We also took advantage of a subsequent suggestion of Kottwitz \cite{KO} to make the group locally compact. The original proposal of Langlands was for a complex proalgebraic group, the category best suited for comparison with Grothendieck's motivic Galois group $\UGal$.

Our main theme ultimately takes the form of an \ubf{explicit} homomorphism
\begin{equation*}
    L_F\longrightarrow\UGal
\end{equation*}
that has two basic interpretations. One is as the general analogue of the Shimura-Taniyama-Weil conjecture, mentioned above. The other is symbolic. As an $L$-homomorphism from the locally compact group $L_F$ to the complex reductive (proalgebraic) group $\UGal$, the map resembles a universal version of the $L$-homomorphisms
\begin{equation*}
    W_F \rightarrow\LG_F
\end{equation*}
from the global Weil group to the $L$-group of a group $G = G_F$ over $F$ that are the objects that Langlands first proposed in \cite{L2} to index some automorphic representations of $G_F$. However, this second interpretation applies only to the automorphic representations that are relevant to arithmetic algebraic geometry.

Some discussion of \cite{Ar2} appeared in Section 9 of the Abel Prize exposition \cite{N1} of the work of Langlands. However, the explicit constructions of \cite{Ar2} are still not widely known. The discussion of them here will be briefer, and somewhat more expository. Our main interest is in the striking implications the constructions might have for the arithmetic geometry and the theory of motives. In this way, the paper here might serve as a supplement to the larger article \cite{N2} in preparation.

Let me add that there are topics in this article that I have yet to fully understand. There may be some misstatements, even though I am confident in the main ideas. It seemed most important to draw attention to these vast areas of mathematics, with their many fundamental questions that seem ripe for study.

\section{Motives}\label{sec: motives}

We begin with a few words on Grothendieck's conjectural theory of motives. Their existence is still dependent on his standard conjectures, despite much progress on them since they were put forward in the $1960s$. However, I think it is fair to say that mathematicians universally regard the conjectures as true, and generally devote their energies to the many extraordinary consequences they will bring.

At its most basic, algebraic geometry is the study of solutions of $m$-polynomial equations in $n$ variables, over a subfield $F\subset\mC$ of the complex numbers. It is arguably the most fundamental form of geometry, whose study goes back centuries.

A simultaneous solution of the equations naturally gives a subset $X$ of $\mC^n$, known as an \ubf{algebraic variety}. Algebraic geometry would then be the study of algebraic varieties, and their subvarieties, known as \ubf{algebraic cycles} on $X$, given as solutions of other sets of polynomial equations in $n$-variables. This is a natural approach to geometry, but there is a basic problem. In general, it is difficult to describe algebraic varieties in explicit terms. For example, there are fundamental ingredients of complex algebraic varieties that algebraic geometry inherits from topology - their cohomology groups, with coefficients in a separate field $Q$? What can we say about these embedded objects for specific algebraic varieties? The answer is, \say{not much}, if we are restricted to the study of their defining equations.

This basic difficulty has what we could regard as a miraculous resolution. Using primarily his extraordinary powers of abstract reasoning, Grothendieck discovered a hidden internal structure to algebraic varieties, which he called \ubf{motifs}, following the use of the term in art and music for unspoken ideas that govern what one sees and hears. Translated as \ubf{motives} in English, these objects can be regarded as the basic ingredients of geometry. When $F$ is a number field and $Q$ is the rational field $\mQ$, motives also assume a fundamental place in number theory, under the term \ubf{arithmetic algebraic geometry}.

As geometric entities, motives can be thought of as fundamental building blocks of algebraic varieties. If we were to think of algebraic varieties as the everyday objects of the world around us, (irreducible) motives could be regarded as the fundamental particles of which they are composed. On the other hand, motives also represent a universal cohomology theory. This became a central question in arithmetic geometry with the introduction of several distinct cohomology theories, such as Betti cohomology from topology, de Rham cohomology from differential geometry (as modified algebraically by Grothendieck for fields $F\subset\mC$), $\ell$-adic \'etale cohomology from arithmetic, Hodge theory from analysis, and so on. The Grothendieck conjectures imply that the motive attached to a variety (as a sum of irreducible motives, together with its decomposition acquired from the underlying weights) represents a universal cohomology theory through which all the others factor.

Much of Grothendieck's original restructuring of algebraic geometry can be summarized as a sequence
\begin{equation}\label{eq: 1.1}
    \begin{tikzcd}
        \cP(F) \ar[r] & \mathcal{E}_{\sim}(F) \ar[r] & \fM^{\scalebox{0.7}{\mbox{eff}}}_{\sim}(F)  \ar[r] & \fM_{\sim}(F) \ar[r] & g\mbox{Vec}_Q
    \end{tikzcd}
\end{equation}
of functors. We present this only as background. It is from the middle row of diagram (6.2) described in Section 6 of \cite{N2}, or in the Section 4 of the general introduction \cite{Andre1}. The left hand term $\cP(F)$ denotes the category of smooth projective varieties over a field $F\subset\mC$. These are algebraic varieties whose complex points form compact complex manifolds.\footnote{They can of course be attached to the affine algebraic varieties above by expanding the affine coordinates $(x_1,\dots, x_m)$ to homogeneous projective coordinates $x = (x_0, x_1,\dots, x_m)$ with the underlying assumption that the nonsingularity condition of the implicit function theorem holds at every point $x$.}
On the right $\gvec{Q}$ denotes the category of finite dimensional graded vector spaces, over a field $Q$ which we usually take to be $\mQ$.

The three interior categories in \eqref{eq: 1.1} represent the three steps in Grothendieck's conjectural transformation of varieties $X\in\cP(F)$ to motives. They are governed by the groups $\mathfrak{Z}^{*}(X)$ of algebraic cycles on varieties $X$, graded by codimension, and their resulting commutative $Q$-algebras
\begin{equation*}
    \mathfrak{Z}^{*}_{\sim}(X) = \mathfrak{Z}^{*}_{\sim}(X)_Q = \mathfrak{Z}^{*}(X)_{Q}/\sim = \mathfrak{Z}^{*}(X)\otimes_{\mZ} Q/\sim,
\end{equation*}
defined by their intersection product and what is called an \ubf{adequate equivalence} relation $\sim$ (\cite[Section 3.1]{Andre1}). The three functors on the left of \eqref{eq: 1.1} are called \ubf{linearization}, \ubf{pseudo-abelianization} and \ubf{inversion}. The first one supplements \ubf{morphisms} $\mbox{Hom}(X, Y)$ in $\cP(F)$ by adjoining elements in the algebra $\mathfrak{Z}^{*}_{\sim}(X \times Y)$ of degree $0$. It takes $\cP(F)$ to an $Q$-linear category $\cE_{\sim}(F)_Q$, which represents the \say{enumerative geometry} practiced by earlier algebraic geometers \cite[Section 1.1]{Andre1}. The second one supplements the objects $X$ in $\cE_{\sim}(F)$ by idempotent elements in $\mathfrak{Z}^{\dim X}_{\sim}(X\times X)$, which is to say, elements in this algebra such that $e^2 = e$. It sends $\cE_{\sim}(F)$ to a pseudo-abelian category $\fM^{\scalebox{0.7}{\mbox{eff}}}_{\sim}(F)$, called the category of \ubf{effective motives}. The third functor supplements the effective motives by adding negative \ubf{weights}. That is, it enlarges $\fM^{\scalebox{0.7}{\mbox{eff}}}_{\sim}(F)$ by expanding the weights it inherits from the cohomology of varieties in $\cP(F)$ in a compatible way from $\mZ^{\ge 0}$ to $\mZ$. It sends $\fM^{\scalebox{0.7}{\mbox{eff}}}_{\sim}(F)$ to the rigid $\otimes$-category $\fM_{\sim}(F)$, known simply as the category of motives.

The category $\fM_{\sim}(F)$ on the right of \eqref{eq: 1.1} thus becomes the object of special interest. However, it still depends on a choice of \say{adequate} equivalence relation \say{$\sim$} on algebraic cycles. The two such relations of particular interest are \ubf{rational equivalence} (the finest) and numerical equivalence (the least fine). The corresponding categories
\begin{equation*}
     \mbox{ \mbox{CHM}}(F) =  \fM_{\scalebox{0.7}{\mbox{rat}}}(F),  
\end{equation*}
and
\begin{equation*}
    \mbox{ \mbox{NM}}(F) = \fM_{\scalebox{0.7}{\mbox{num}}}(F), 
\end{equation*}
called \ubf{Chow motives} and \ubf{numerical motives} respectively, then come with a functor
\begin{equation}\label{eq: 1.2 FC}
      \mbox{CHM}(F)\longrightarrow  \mbox{NM}(F).
 \end{equation}
This completes the preliminary discussion of motives as building blocks. They take the form of the Tannakian category of numerical motives $\mbox{ \mbox{NM}}(F)$.

However, we are equally concerned with how Grothendieck relates motives to various cohomology theories from geometry. As in the original volume \cite{R} from 1972, one can define a Weil cohomology axiomatically as a functor
\begin{equation}\label{eq: FC 1.3}
     H^*: \mathcal{P}(F)^{\scalebox{0.7}{\mbox{op}}}\longrightarrow \mbox{Vect}_L^{\ge 0},
\end{equation}
that satisfies conditions $1)$ and $2)$ in \cite[3.1.1.1]{Andre1} (as well as the basic normalization 
\begin{equation*}
    H^i(\mathbb{P}^1) \cong \begin{cases}
        L, & \mbox{ if } i = 0, 2,\\
        0, & \mbox{ otherwise ,}
    \end{cases}
\end{equation*}
for projective $1$-space $\mathbb{P}^1$). Here $K$ is any field containing $F$, and $\mathcal{P}(F)^{\scalebox{0.7}{\mbox{op}}}$ is the opposite of the category $\cP(F)$, which reverses arrows in morphisms to accommodate cohomology as a contravariant functor. In \cite[4.2.1-4.2.5]{Andre1}, this definition is compared with the contravariant functor
\begin{equation*}
    \mathfrak{h} = \mathfrak{h}_{\sim} : \mathcal{P}(F)^{\scalebox{0.7}{\mbox{op}}}\longrightarrow \fM_{\sim}(F)
\end{equation*}
obtained from the composition of the first three functors in \eqref{eq: 1.1}. The conclusion is that the choice of a Weil cohomology is equivalent to that of a $\otimes$-functor 
\begin{equation}\label{eq: FC 1.4}
    H^*: \mbox{CHM}(F)_Q\longrightarrow (g\mbox{Vec})_L, 
\end{equation}
such that $H^i(\mathbb{P}^1) = 0$ for $i\neq 2$ \cite[Proposition 4.2.5.1]{Andre1}. The argument also implies that \eqref{eq: FC 1.4} factors to a (unique) composition
\begin{equation}\label{eq: FC 1.5}
    \mbox{CHM}(F)_Q\longrightarrow \mbox{NM}(F)_Q \longrightarrow (g\mbox{Vec})_L,
\end{equation}
where the first arrow is the projection \eqref{eq: 1.2 FC}. 

We can therefore think of $ \mbox{CHM}(F) = \mbox{CHM}(F)_Q$ as a universal cohomology for smooth projective varieties. It is not a Weil cohomology itself, but rather a category that classifies them through their realizations \eqref{eq: FC 1.4}. Nor do its objects represent the motives we would like to work with. This distinction is reserved for the Tannakian category of numerical motives, in which the right hand arrows in \eqref{eq: 1.1} and \eqref{eq: FC 1.5} represent a fibre functor in the form of Betti cohomology. We should perhaps add that $g\mbox{Vec}_L$ is the graded category attached to the Tannakian category $\mbox{Vec}_L$ of finite dimensional vector spaces over $L$. It is a general operation, which modifies the commutativity of tensor products by certain signs (Koszul signs) for objects of negative degree \cite[Example 2.2.2.2]{Andre1}, and has the effect of making the compositions in \eqref{eq: 1.1} and \eqref{eq: FC 1.5} \ubf{contravariant functors}. Chow motives can thus be regarded as refinements of motives. The fibres in $\mbox{CHM}(F)$ from \eqref{eq: FC 1.5} provide the cohomology of motives $M$ in $\mbox{NM}(F)$, plus interesting supplementary information. With this observation, we have completed our preliminary discussion of cohomology.

\section{Automorphic representations}\label{sec: automorphic representations}

Automorphic representations are the foundation of the Langlands program. They were introduced for general groups $G$ in 1959 by Harish-Chandra \cite{N3}, \cite{N4} in the closely related setting of automorphic forms. Special cases go back further to Poincar\'e, Hilbert-Blumenthal, Siegel, Maas, Selberg and others, especially for the group $\SL(2)$, where they are widely known as modular forms. But the present day theory was introduced by Langlands. He presented it in three distinct formats: his 1965 volume on the spectral classification of general Eisenstein series \cite{N5}, his 1967 letter to Weil \cite{N6}, and his revolutionary 1969 lecture \cite{L2} on seven \say{Problems in the theory of automorphic forms}.

Since our goal is to describe motives and automorphic representations as theories that at the deepest levels are intertwined, let us introduce a diagram
\begin{equation}\label{eq: 2.1}
    \begin{tikzcd}
        \{G\} \ar[r] & \{\pi\} \ar[r] & ^LG = \widehat{G}\rtimes \Gamma_{F} \ar[r] & (r, V_r) \ar[r] & L(s, \pi, r)
    \end{tikzcd}
\end{equation}
that can play the role of \eqref{eq: 1.1}. It is more impressionistic, to be sure, and its ties to the terms in \eqref{eq: 1.1} are not to be taken literally, but it does give us some perspective on the theory of Langlands.

On the left, $G$ ranges over (connected) reductive algebraic groups defined over a fixed number field $F$. These have a reasonable classification defined by Coxeter - Dynkin diagrams, and twists by Galois $1$-cocycles that can be resolved in terms of abelian class field theory. The second term consists of (irreducible) automorphic representations of the adelic group $G(\mA)$. They are the primary objects of study. The third term is the Langlands $L$-group, defined as a semidirect product
\begin{equation*}
    ^LG = \widehat{G}\rtimes \Gamma_{F},\mkern70mu \Gamma_F = \mbox{Gal}(\overline{F}/F),
\end{equation*}
of the complex connected dual group $\whG$ of $G$ and the absolute Galois group of $F$. It was Langlands' expectation, expressed in \cite{N6} and \cite{L2}, that the complex coordinates in $\whG$ would represent parameters for representations. The fourth term is a finite dimensional representation
\begin{equation*}
    r:\LG \longrightarrow\GL(V_r).
\end{equation*}
This transfers data for $\pi$ from $\LG$ to the complex general linear group $\GL(V_r)$. The fifth term is an automorphic $L$-function
\begin{equation}\label{eq: 2.2}
    L(s, \pi, r) = \prod_v L_v(s, \pi_v, r), \SP  s>>1,
\end{equation}
attached to $\pi$ and $r$. He defined it as an absolutely convergent Euler product, the multiplication of whose factors becomes a completely new Dirichlet $L$-series. He conjectured in \cite{N6} and \cite{L2} that $L(s, \pi, r)$ has meromorphic continuation in $s$ to the complex plane, with the kind of functional equation familiar from number theory. He suggested that the family of $L$ functions so obtained was universal, in the sense that it includes all arithmetic $L$-functions that arise from number theory and arithmetic algebraic geometry.

As in the case \eqref{eq: 1.1} for motives, we will comment only briefly for background on the terms of the sequence \eqref{eq: 2.1} for automorphic representations. We recall that if $G$ is a reductive algebraic group over a number field $F$, its group of $F$-rational points comes with a natural embedding
\begin{equation*}
    G(F) \subset G(\mA),
\end{equation*}
as a discrete subgroup of the locally compact group $G(\mA)$ of adelic points. This stems from the fact that
\begin{equation*}
    \mA = \mA_F = \displaystyle\widetilde{\bigotimes_v}\, F_v
\end{equation*}
is a locally compact ring, defined as a restricted direct sum of all the local completions of $F$, in which the diagonal embedding of $F$ is a discrete, co-compact subring. Adelic groups were not widely used when Langlands announced his conjectures, but they were eventually seen to provide the most efficient setting. The point is that there are data embedded in arithmetic quotients of real groups, such as Hecke correspondences, that become more transparent and easier to work with in adelic quotients.

In his original paper \cite{L2}, Langlands introduced what later became known as an \ubf{automorphic} representation of $G(\mA)$ (or just $G$). He defined it informally as an irreducible representation of the group $G(\mA)$ that \say{occurs in} the decomposition of $L^2(G(F)\backslash G(\mA))$ under right translation by $G(\mA)$. Stated thus, the definition is understood implicitly in terms of the decomposition of a representation of $G(\mA)$ on a Hilbert space into irreducible unitary representations. It is imprecise for two reasons. One is that $L^2(G(F)\backslash G(\mA))$ generally has a part that decomposes continuously (like the theory of Fourier transforms) whose irreducible constituents are defined only up to a set of measure $0$. The other is that automorphic representations in the continuous part can have analytic continuation in their continuous parameters to nonunitary representations, and these should also be regarded as automorphic. The formal definition of an automorphic representation came later, in the Corvallis proceedings \cite{BJ} and \cite{L5}.

In any case, the informal definition was quite sufficient for the new phenomena Langlands was about to describe. An automorphic representation $\pi\in\Pi_{aut}(G)$ would come with a decomposition into a restricted tensor product
\begin{equation*}
    \pi = \displaystyle\widetilde{\bigotimes_v}\, \pi_v
\end{equation*}
of irreducible representations $\pi_v\in \Pi(G_v)$ of the local completions $G_v = G(F_v)$, almost all of which
\begin{equation}\label{eq: 2.3}
    \{\pi_v : v\notin S\}, \SP S \mbox{ finite }\supset S_{\infty},
\end{equation}
are \ubf{unramified}. Unramified representations $\{\pi_v\}$ are a particularly simple class of irreducible representations of the local groups $\{G_v\}$. Without recalling the precise definition, we can say that they are attached to places $v$ for which $G_v$ is \ubf{unramified} (in a natural sense), and that they are then irreducible constituents of representations induced from unramified quasi-characters $\chi_v$ on a Borel subgroup $B_v(F_v)$ of $G_v(F_v)$. Langlands' $L$-group is defined so as to map these quasi-characters to semisimple elements in $\LG$. The unramified representations $\pi_v$ of $G_v$ are then bijective with semisimple conjugacy classes $c(\pi_v)$ in $\LG$ that are \ubf{unramified}, in the sense that their image in $\Gamma_F$ is the Frobenius conjugacy class at $v$. An automorphic representation $\pi$ thus provides a concrete set of data
\begin{equation}\label{eq:2.4}
    c = c(\pi) = \{c_v = c(\pi_v): v\notin S\},
\end{equation}
where $S$ is a finite set of valuations of $F$ that contains the set $S_{\infty}$ of archimedean valuations, and $c(\pi_v)$ is an unramified conjugacy class in$\LG$.

With this understanding, the automorphic $L$-function introduced by Langlands was very natural. For any $r$ and $\pi$, he formulated it as an Euler product
\begin{align*}
    L(s, \pi, r) &= \prod_v L_v(s, \pi, r) = \prod_{v\in S} L(s, \pi_v, r_v) \;\cdot\; \prod_{v\notin S} L(s, \pi_v, r_v)\\
    &= L_S(s, \pi, r)\cdot L^S(s, \pi, r),
\end{align*}
in which for almost all valuations $v\notin S$ of $F$, the local factors were defined
\begin{equation*}
    L(s, \pi_v, r_v) = \det(1 - r_v(c_v)\cdot q_v^{-s})^{-1}
\end{equation*}
in terms of the characteristic polynomials of corresponding conjugacy classes
\begin{equation*}
    r_v(c_v) = r(c(\pi_v)).
\end{equation*}
This is a natural generalization of both Dirichlet and Artin $L$-functions, where $q_v$ is as usual the order of the residue field of the $v$-adic completion $F_v$ of $F$.

Langlands' foundational initial question from \cite{L2} can then be stated as follows.
\begin{description}[leftmargin = 0pt]
    \item[Langlands Question 1 \cite{L2}:] \textit{Is it possible to define local $L-$functions
    \begin{equation*}
        L_v(s, \pi, r) = L(s, \pi_v, r_v)
    \end{equation*}
    and local $\varepsilon$-factors
    \begin{equation*}
        \varepsilon_v(s, \pi, r, \psi) = \varepsilon(s, \pi_v, r_v, \psi_v) = \varepsilon(\pi_v, r_v, \psi_v)q_v^{-(n_v - \frac{1}{2}s)},
    \end{equation*}
    for valuations $v\in S$ and any nontrivial additive character $\psi$ on $F\backslash\mA$, such that the Euler product $L(s, \pi, r)$ has analytic continuation, and satisfies the functional equation
    \begin{equation*}
        L(s, \pi, r) = \varepsilon(s, \pi, r)L(1 - s, \pi, r^{\vee}),
    \end{equation*}
    for
    \begin{equation*}
        \varepsilon(s, \pi, r) = \prod_{v\in S} \varepsilon_v(s, \pi, r, \psi)\,\mbox{?}
    \end{equation*}}
\end{description}
\section{Automorphic Galois group}\label{sec: automorphic galois groups}

In this section, we describe our conjectural candidate for the automorphic Galois group $L_F$. Conjecture enters not just in the presumed role of $L_F$, but also in its basic construction. In fact, the definition we give depends directly on Langlands' Principle of Functoriality, a conjecture that is at the heart of the Langlands program, and that was implicit (without the name) throughout the formative paper \cite{L2}. In its most basic global form, it can be stated as follows.

\begin{description}[leftmargin=0pt]
    \item[Langlands Principle of Functoriality:] \textit{Suppose that $G$ and $G'$ are quasisplit groups over the number field $F$, and that
    \begin{equation*}
        \rho: \LG'\longrightarrow \LG
    \end{equation*}
    is an $L$-homomorphism\footnote{That is, a homomorphism that is complex analytic on $\whG'$, and that commutes with the projections of $\LG'$ and $\LG$ onto $\Gamma_F$.}
    between their $L$-groups. Then if $\pi'\in \Pi_{aut}(G')$ is an automorphic representation of $G'$, there exists an automorphic representation $\pi = \rho'(\pi')$ of $G$ such that
    \begin{equation*}
        c(\pi_v) = \rho'(c(\pi_v')), \SP v\notin S.
    \end{equation*}
    In other words, the concrete data in $\pi'$ and $\pi$ given by families of semisimple conjugacy classes in their $L$-groups match under the transfer defined by $\rho'$.}
\end{description}

For an illustration of the power of functoriality, suppose that $L(s, \pi, r)$ is a general Langlands $L$-function. We then apply the assertion of functoriality with $(G, \pi, r, \GL(N_r)), N_r = \dim r$, in place of the quadruple $(G', \pi', G, \rho')$ above. The assertion is that there is an automorphic representation $\pi_r$ of $\GL(V_r)$ such that
\begin{equation*}
    c(\pi_{r, v}) = r(c(\pi_v))
\end{equation*}
for almost all $v\notin S$. This gives an identity
\begin{equation}\label{eq: 3.1}
    L^S(s, \pi, r) = L^{S}(s, \pi_r, S_r)
\end{equation}
of unramified $L$-functions, where $S_r$ is given by the standard, $n$ dimensional representation of $\widehat{\GL}(N_r) = \GL(V_r)$.

The point here is that the conjecture of Langlands' Question 1 already has an affirmative answer for the standard $L$-functions of general linear groups $\GL(N)$. The solution in fact has a long history. It was first established by Hecke for many cases with $N = 1$ \cite{N7} and some cases for $N = 2$ \cite{N8}. It was then established by Tate in general for $N = 1$ \cite{T}, by Jacquet and Langlands in general for $N = 2$ \cite{N9}, and finally by Godement-Jacquet \cite{N10} for any $N$ . In all these cases, the ideas for the proof were roughly similar. In the papers \cite{T} and \cite{N10}, for example, one studies additive Fourier transforms of test functions
\begin{equation*}
    \phi_s(x) = \phi(x)|\det x|^s, \SP \phi\in C_c^{\infty}(M(N, \mA)), x\in M(N, \mA),
\end{equation*}
on adelic matrix space $M(N, \mA)$, with a supplementary complex variable $s\in\mC$. Tate's elegant thesis \cite{T} is particularly clear, and has traditionally been where students begin the study of such things. Then, given a proof of the conjecture for the standard $L$-functions on the right hand side of \eqref{eq: 3.1}, one obtains a proof for the general $L$-functions on the left hand side of \eqref{eq: 3.1} by defining the ramified local $L$-functions $L_v(s, \pi, r)$ and $\varepsilon$-factors $\varepsilon_v(s, \pi, r, \psi)$ from their analogues on the right hand side.

We turn now to the automorphic Galois group $L_F$. Automorphic representations should be governed, perhaps even classified, by certain locally compact groups attached to $F$ and its local completions $F_v$. These groups would fit into a commutative diagram 
\begin{equation}\label{eq: 3.2}
    \begin{tikzcd}
      L_{F_v} \ar[r, mapsto]\ar[d, hook] & W_{F_v} \ar[r, mapsto]\ar[d, hook] &\Gamma_{F_v} \ar[d, hook]\\
      L_{F} \ar[r, mapsto] & W_{F} \ar[r, mapsto] & \Gamma_{F} .
  \end{tikzcd} 
\end{equation}
The vertical arrows represent embeddings defined up to conjugacy. On the right, we have the absolute Galois group $\Gamma_F$ with its embedded localizations $\Gamma_{F_v}$. In the middle of \eqref{eq: 3.2}, we have the absolute global Weil group $W_F$, with embeddings of its local analogues $W_{F_v}$. These objects had been introduced by Weil \cite{N11} in 1951 as enrichments of Galois groups provided by abelian class field theory. But they assumed new importance with Langlands' Bochner lecture \cite{L2}.

Langlands conjectured that there would be automorphic representations $\pi$, for a group $G$ over $F$, attached to $L$-homomorphisms from $W_F$ to its $L$-group$\LG$. This global correspondence would be compatible with the local embeddings in \eqref{eq: 3.2}. For it would give the local constituents $\pi_v$ of $\pi$ in terms of a separate local correspondence for each $v$, which would take $L$-homomorphisms from $W_{F_v}$ to $\LG_v$ to irreducible representations $\pi_v\in\Pi(G_v)$ of $G_v = G(F_v)$. The local correspondence would in fact amount to a classification, conjectured later by Langlands, in which the local Weil group $W_{F_v}$ is replaced by its variant
\begin{equation*}
    L_{F_v} = 
    \begin{cases} 
          W_{F_v}, & \mbox{ if } v\in S_{\infty}, \\
          W_{F_v} \times \mbox{SU}(2), & \mbox{ if } v\notin S_{\infty},
       \end{cases}
\end{equation*}
the Weil-Deligne group in the upper left hand corner of diagram \eqref{eq: 3.2}. Suppose for example that $G_v$ is quasi-split. Then there should be a bijection, from equivalence classes of $L$-homomorphisms
\begin{equation}\label{eq: 3.3}
    \phi_v: L_{F_v} \longrightarrow \LG_v,
\end{equation}
to finite disjoint packets $\Pi_{\phi_v}$, whose union is the full set $\Pi(G_v)$. It would come with further structure that would characterize the finite packets in terms of Langlands' local conjecture of endoscopy \cite{N12}. This \say{local Langlands' correspondence} was established for real groups $G_v, v\in S_{\infty}$, in \cite{N13}. It has been established more recently for quasisplit $p$-adic classical groups $G_v$ \cite{Ar3}, and then unitary groups \cite{Mo}.

However, the conjectured global Langlands correspondence, for maps
\begin{equation*}
    \phi: W_F\longrightarrow\LG,
\end{equation*}
leaves out many automorphic representations. The problem is that the global Weil group is too small. To complete it, one would have to replace $W_F$ by the general enriched group $L_F$ in the lower left hand corner of the diagram.

The construction of $L_F$ is based on the subset of \ubf{primitive} automorphic representations attached to certain groups $G$ over $F$. Taking this notion for granted for the moment, we identify the two ingredients that go into the construction of $L_F$ as follows.

(i). An indexing set
\begin{equation*}
    \cC_F = \{(G, c)\},
\end{equation*}
where $G/F$ is a quasisplit, simple, simply connected group over $F$, and
\begin{equation*}
    c = c(\pi) = \{c_v = c(\pi_v) \;:\; v\notin S\}
\end{equation*}
is an equivalence class of data \eqref{eq:2.4} attached to a primitive automorphic representation $\pi$ of $G$ (where two data $c = \{c_v\}$ and $c' = \{c_v'\}$ are \ubf{equivalent} if there is an isomorphism from $G$ to $G'$ over $F$ whose dual takes $c_v$ to $c_v'$ for almost all $v$. See the definition on page 472 of \cite{Ar2}).

(ii). For each $c\sim (G, c)$ in $\cC_F$, an extension
\begin{equation}\label{eq: 3.4}
    1 \longrightarrow K_c \longrightarrow L_c \longrightarrow W_F \longrightarrow 1
\end{equation}
of $W_F$ by a compact real form $K_c$ of the simply connected cover $\whG_{sc}$ of the complex dual group $\whG$. (See \cite{Ar2}, where the extension is stated as (3.3), and constructed in Section 4).

Given (i) and (ii), we define the automorphic Galois group simply as the locally compact fibre product
\begin{equation}\label{eq: 3.5}
    L_F = \prod_{c\in\mathcal{C}_F}(L_c\longrightarrow W_F),
\end{equation}
over $F$. It is an extension
\begin{equation}\label{eq: 3.6}
    1 \longrightarrow K_F \longrightarrow L_F \longrightarrow W_F \longrightarrow 1
\end{equation}
of the locally compact group $W_F$ by the compact group
\begin{equation*}
    K_F = \prod_{c\in\mathcal{C}_F}K_c,
\end{equation*}
and is therefore a locally compact group itself. It should come with a localization
\begin{equation*}
      \begin{tikzcd}
       L_{F_v} \ar[r, mapsto]\ar[d, hook] & W_{F_v} \ar[d, hook]\\
       L_c \ar[r, mapsto] & W_F
   \end{tikzcd}
\end{equation*}
at any valuation $v$ of $F$ for each $c$ in the product \eqref{eq: 3.5} for $L_F$, and hence a localization that gives the left hand side square in \eqref{eq: 3.2}.

Our construction is thus founded on a notion of primitive automorphic representations $\pi\in\Pi_{prim}(G)$. In describing it, our first point is that $G$ is restricted as above to being a quasisplit, simple, simply connected group over $F$. We should note here that $G$ is not required to be \ubf{absolutely} simple. It is thus typically of the form $\Res_{E/F}(G_E)$, where $G_E$ is an absolutely simple group over some extension $E/F$. With this understanding, we begin to see the complexity of the group $L_F$. In particular, its internal structure is closely tied to general base change of automorphic representations of $G$. The next point is that for the given $G$, there are actually two characterizations of primitive representations that turn out to be equivalent. The first is that $\pi$ is not a proper functorial image, which is to say that it is not the functorial image of a pair $(G', \pi')$ with $\dim G' < \dim G$. The other is that $\pi$ is \ubf{tempered}, which is to say that its local constituents $\pi_v$ all come from local parameters $\phi_v\in\Phi(G_v)$ with \ubf{bounded} image, and in addition, that for any $r$, the unramified $L$-function
\begin{equation*}
    L(s, c^S, r) = L^S(s, \pi, r), \SP c^S = c^S(\pi),
\end{equation*}
has minimal order at $s = 1$. In other words
\begin{equation*}
    \mbox{ord}_{s = 1}(L(s, c^S, r)) = [r : \mathbf{1}_{^LG}],
\end{equation*}
where $\mathbf{1}_{^LG}$ is the trivial one dimensional representation of$\LG$.

Given now the conjectural construction of $L_F$, we can state its main conjectural property.
\begin{conjecturegroup}\label{conj: 1}

    \begin{description}[leftmargin=0pt]
        \item[(a).] There is a bijection
         \begin{equation*}
        \{r: L_F\longrightarrow \mbox{GL}(N, \mathbb{C})\} \xlongrightarrow{\sim} \{\pi\in \Pi_{cusp, temp}(\mbox{GL}(N))\},
    \end{equation*}
    from equivalence classes of bounded, irreducible, $N$-dimensional representations of $L_F$ to cuspidal tempered automorphic representations of $\GL(N)$, which is compatible with the local Langlands correspondence at each $v$.
        \item[    (b).] For any quasisplit group $G$ over $F$, there is a bijection
        \begin{equation*}
        \{\phi: L_F\longrightarrow\,^LG\} \xlongrightarrow{\sim} \{\Pi_{\phi}\subset \Pi_{aut,\, temp}(G)\},
        \end{equation*}
    from equivalence classes of bounded global parameters $\Phi_{bdd}(G)$ to disjoint global packets, whose union contains the set $\Phi_{aut, temp}(G)$ of tempered automorphic representations of $G(\mA)$, and which is compatible with the local Langlands correspondence at each $v$.
    \end{description}
\end{conjecturegroup}
\begin{description}[leftmargin=0pt]
        \item[Remarks 1.] The global form of Langlands' conjectural theory of endoscopy comes with an explicit formula for the multiplicity of any representation $\pi\in\Pi_{\phi}$ in the tempered automorphic spectrum of $G(\mA)$. This characterizes the automorphic image of the conjectural global Langlands correspondence.
        \item[2.] The mappings (a) and (b) do not give all representations in the automorphic discrete spectrum. They do not in general even account for all cuspidal automorphic representations, since for some $G$ these are not always tempered. To account for the full automorphic discrete spectrum, one must enlarge $L_F$ to a product $A_F = L_F\times \SL(2, \mC)$, and work with $A$-packets $\Pi_{\psi}$ for parameters
        \begin{equation*}
            \psi: A_F\longrightarrow\LG
        \end{equation*}
        in place of the global $L$-packets $\Pi_{\phi}$ above.
\end{description}

\section{Motivic Galois group}\label{sec:motivic galois groups}

There are qualitative differences between the motivic Galois group $\UGal$ and the automorphic Galois group $L_F$. For a start, $\UGal$ is a complex proalgebraic group rather than a locally compact topological group. Its complex simple factors are thus larger than their compact real forms that give corresponding factors for $L_F$. On the other hand $\UGal$ is smaller in that $L_F$ has compact simple factors that do not contribute to $\UGal$.

Despite these differences, our construction of $\UGal$ is parallel to that of $L_F$. It takes the form of a short exact sequence that embeds into a commutative diagram 
\begin{equation}\label{eq: 4.1}
    \begin{tikzcd}
        1 \ar[r] & K_F \ar[r] \ar[d] & L_F \ar[r]\ar[d] & W_F \ar[r]\ar[d] & 1\\
        1 \ar[r] & \cD_F \ar[r] & \UGal \ar[r] & \cT_F \ar[r] & 1
    \end{tikzcd}
\end{equation}
in which the top row is the extension \eqref{eq: 3.6} that gives $L_F$, and the bottom row is its analogue for $\UGal$. We are following \cite[Section 3]{N2} here, which is based on the original article \cite{Ar2}. The complex motivic Galois group is at the centre of the bottom row of \eqref{eq: 4.1}. The three vertical arrows represent continuous mappings from locally compact groups to complex proalgebraic groups. We shall describe them.

Recall that $L_F$ is a locally compact fibre product over a set $\cC_F$ of indices $(G, c)$. The group $\UGal$ will be a proalgebraic fibre product
\begin{equation}\label{eq: 4.2}
    \prod_{c\in\mathcal{C}_{F,0}}(\mathcal{G}_c\longrightarrow \mathcal{T}_F)
\end{equation}
for a subset $\cC_{F, 0}$ of indices in $\cC_F$ we designate to be of type $A_0$, and for each index $c\sim (G, c)$ in $\cC_{F, 0}$, an extension of proalgebraic groups
\begin{equation}\label{eq: 4.3}
    1 \longrightarrow \mathcal{D}_c \longrightarrow \mathcal{G}_c \longrightarrow \mathcal{T}_F \longrightarrow 1.
\end{equation}
We shall describe these two ingredients in turn, requiring somewhat more detail than we needed for their $L_F$ analogues in the last section.

(i). Suppose that
\begin{equation*}
    c\sim (G, c), \SP c = c(\pi) = \{c_v = c(\pi_v) : v\notin S\},\;\; \pi\in\Pi_{prim}(G),
\end{equation*}
represents an index in $\cC_F$. The condition that $c$ be of type $A_0$ is given in terms of the local Langlands parameters $\phi_v$ for the archimedean components $\pi_v$ of $\pi$, which is to say the $L$-homomorphisms $\phi_v$ in \eqref{eq: 3.3} for $v\in S_{\infty}$ such that $\pi_v$ lies in the local $L$-packet $\Pi_{\phi_v}$. At first glance this does not seem plausible. Indeed, for the index $c$, $\pi$ is any primitive automorphic representation such that $c_v$ equals $\pi(c_v)$ for all $v\notin S$. How do we know that $c$ then determines the archimedean components of $\pi$, or at least the Langlands parameters
\begin{equation*}
    \{\phi_v : v\in S_{\infty}\}
\end{equation*}
for these components? The answer is an implicit conjecture. For any primitive automorphic representation $\pi$ of the group $G$, we conjecture that \ubf{all} of its local Langlands parameters are uniquely determined by its unramified parameters. In particular, the data in the index $c$ determines a set of archimedean parameters as above. This would be an analogue for primitive automorphic representations of $G$ of the theorem of strong multiplicity $1$ for cuspidal automorphic representations of general linear groups.

Taking this conjecture for granted, we define $\cC_{F, 0}$ to be the set of indices in $\cC_F$ that are of Hodge type for each $v\in S_{\infty}$. In other words, for each $v\in S_{\infty}$ and every finite dimensional (continuous) representation $r$ in $\LG$, the associated representations $r\circ\phi_v$ of archimedean Weil groups $W_{F_v}$ are of Hodge type, in the sense that the restriction of $r\circ\phi_v$ to the subgroup $\mC^*$ of $W_{F_v}$ is a direct sum of quasicharacters of the form
\begin{equation*}
    z \longrightarrow z^{-p}\overline{z}^{-q}, \SP z\in\mathbb{C^*}, \;p, q\in\mathbb{Z}.
\end{equation*}
(ii). Our second ingredient is the short exact sequence \eqref{eq: 4.3}, for an index $c\sim (G, c)$ now in $\cC_{F, 0}$. The group $\cT_F$ on the right of \eqref{eq: 4.3} is Langlands' Taniyama group, defined in Section 5 of his Corvallis proceedings \cite{L4} as an explicit extension
\begin{equation}\label{eq: 4.4}
    1 \longrightarrow \mathcal{S}_F \longrightarrow \mathcal{T}_F \longrightarrow \Gamma_F \longrightarrow 1
\end{equation}of proalgebraic groups. For this, $\cS_F$ is the Serre group over $F$, a complex proalgebraic torus defined over $\mQ$ \cite{Se1}, \cite{Se2}, \cite{L4}. As we would expect, $\cS_F$ is determined by its lattice $X_F$ of rational characters, described in explicit terms in \cite[Section 5]{L4} and \cite[Section 3]{Se2}. Given $\cS_F$, the Taniyama group $\cT_F$ is then defined by an explicit, but quite elaborate Galois $2$-cocycle on $\Gamma_F$, with values in $\cS_F$. We can think of $\cT_F$ intuitively as the \say{algebraic hull} of the \say{motivic part} of the Weil group $W_F$. We illustrate this by embedding the extension \eqref{eq: 4.4} into a commutative diagram
\begin{equation}\label{eq: 4.5}
    \begin{tikzcd}
    1 \ar[r] & W_F^0 \ar[r] \ar[d] & W_F \ar[r]\ar[d] & \Gamma_F \ar[r] \ar[d, phantom, "||"]& 1\\
    1 \ar[r] & \mathcal{S}_F \ar[r] & \cT_F \ar[r] & \Gamma_F \ar[r] & 1       
    \end{tikzcd}
\end{equation}
of short exact sequences, where $W_F^0$ is the connected component of $1$ in the global Weil group.

The group $\cD_c$ on the left of \eqref{eq: 4.3} is the complexification of the compact group $K_{c}$ in the corresponding extension for $L_c$. In more direct terms, it is just the simply connected covering group $\whG_{sc}$ of the dual group $\whG$ of $G$. The groups $\cT_F$ and $\cD_c$ are the same as in \cite[Section 6]{N2}. The required extension \eqref{eq: 4.3} is then just \cite[(6.3)]{N2}, which is given by the associated definition \cite[(6.4)]{N2} of the group $\cD_c$. In conclusion, we observe that the definition \cite[(6.4)]{N2} of the algebraic extension \eqref{eq: 4.3} for the algebraic group $\mathcal{G}_c$ is completely parallel to the definition \cite[(4.4)]{N2} of our extension \eqref{eq: 3.3} for the locally compact group $L_c$.

Given (i) and (ii), we define the complex motivic Galois group as the fibre product \eqref{eq: 4.2}. It is an extension
\begin{equation}\label{eq: 4.6}
    1 \longrightarrow \mathcal{D}_F \longrightarrow \mathcal{G}_F \longrightarrow \mathcal{T}_F \longrightarrow 1
\end{equation}
of $\cT_F$ by the proalgebraic group
\begin{equation}\label{eq: 4.7}
    \mathcal{D}_F = \prod_{c\in\mathcal{C}_{F,0}}\mathcal{D}_c.
\end{equation}
The extension \eqref{eq: 4.6} can then take its place as bottom row of the commutative diagram \eqref{eq: 4.1}. Notice the similarity of \eqref{eq: 4.1} with the commutative diagram \eqref{eq: 4.5} for $\cT_F$. Extending the remark preceding \eqref{eq: 4.5}, we can think of $\UGal$ intuitively as the \say{algebraic hull} of the \say{motivic part} of $L_F$. In fact, \eqref{eq: 4.5} represents a part of this general process, namely its restriction to the complement of the connected, semisimple subgroup of the disconnected, proreductive group $\UGal$.

Another comment concerns the mapping from $L_F$ to $\UGal$ at the centre of the diagram \eqref{eq: 4.1}. This is to be regarded as the general analogue of the Shimura-Taniyama-Weil conjecture mentioned in the foreword. We can imagine that its proof, whenever it comes, will be tied up with the ultimate proof of Langlands functoriality.
\section{$\mQ$-structure}\label{sec: Q structure}

We have just described the group $\UGal$ explicitly as a complex, reductive, disconnected, proalgebraic group, which depends on a number field $F$ with a complex embedding $F\subset\mC$. Grothendieck's conjectures require that it be defined over $\mQ$. Can one make this explicit? A related question is how it depends on the underlying number field $F\subset \mC$. The answer to the first question is undoubtedly yes, but there are some details that have still to be resolved. Some of these are in the forthcoming paper \cite{N2}, while others will appear all being well, in a paper \cite{N14} under construction.

The second question applies to $L_F$ as well as $\UGal$. It is in fact natural to describe it in this case first. For $L_F$ is a locally compact group, which was defined in Section \ref{sec: automorphic galois groups} as an extension \eqref{eq: 3.5} of the global Weil group $W_F$. The dependence of $W_F$ on $F$ is fundamental to its existence. It is described on the first pages of Tate's Corvallis article on the number theoretic background for automorphic representations. The property in question is a bijection 
\begin{equation*}
    \begin{tikzcd}
        W_{\mQ}/W_F \arrow{r}{\sim} & \Gamma_{\mQ}/\Gamma_F \ar[r, "\sim"] & \mbox{Hom}(F, \mC)
    \end{tikzcd}
\end{equation*}
that reduces the dependence of $W_F$ on $F$ to the same property for Galois Groups. It is natural to conjecture that this property extends to $L_F$.

\refstepcounter{conjecturegroup}
\setcounter{subconjecture}{0}

\begin{subconjecture}\label{conj: 2}
    The automorphic Galois groups $L_F$ and $L_{\mQ}$ fit into a commutative diagram of continuous homomorphisms
        \begin{equation*}
    \begin{tikzcd}
    L_F \ar[r, mapsto]\ar[d, hook] & W_F \ar[r, mapsto]\ar[d, hook] & \Gamma_F\ar[d, hook]\\
     L_{\mQ} \ar[r, mapsto] & W_{\mQ} \ar[r, mapsto] & \Gamma_{\mQ}        
    \end{tikzcd}
\end{equation*}
in which the coset spaces for the three injections have canonical isomorphisms
    \begin{equation}\label{eq: 5.1}
    \begin{tikzcd}
        L_\mQ/L_F \arrow{r}{\sim} & W_\mQ/W_F \arrow{r}{\sim} & \Gamma_\mQ/\Gamma_F \ar[r, "\sim"] & \mbox{Hom}_{\mQ}(F, \mC).
    \end{tikzcd}
\end{equation}
In particular, if $F/\mQ$ is a Galois extension, $L_{\mQ}/L_F$ is isomorphic to the Galois group $\Gamma_{\mQ/F}$.
\end{subconjecture}

This is Conjecture 4.1 of \cite{N2}. We refer the reader to the discussion in \cite{N2} following the statement of the conjecture. In brief terms, Conjecture \ref{conj: 2} would be a fundamental property of automorphic representations that should be closely tied to the conjectural construction of $L_F$ itself. This is not surprising, given that general base change has also been thought to have a fundamental role in any general proof of functoriality. A starting point at this stage would be the recent extension \cite{CR} of cyclic base change for $\GL(N)$ to solvable base change for $\GL(N)$. It would be very interesting to extend the results in \cite{CR} to other groups $G$, especially to the cases of quasisplit classical and unitary groups, where the general theory of local and global endoscopy has been established \cite{Ar3}, \cite{Mo}, \cite{AGIKMS}. We hope to include further discussions of base change in the paper \cite{N14}.

Conjecture \ref{conj: 2} ought to extend to the motivic Galois group $\UGal$. There are actually two cases to consider. The first applies to $\UGal$ as a complex proalgebraic group. Our discussion of its construction from the last section leads to the following variant of Conjecture \ref{conj: 2}.


\begin{subconjecture}\label{conj: 2(b)}
    The complex motivic Galois groups $\UGal$ and $\mathcal{G}_{\mQ}$ fit into a commutative diagram of proalgebraic homomorphisms 
    \begin{equation*}
    \begin{tikzcd}
        \UGal\ar[r, mapsto]\ar[d, hook] & \cT_F \ar[r, mapsto]\ar[d, hook] & \Gamma_F\ar[d, hook]\\
        \mathcal{G}_{\mQ}\ar[r, mapsto] & \cT_{\mQ} \ar[r, mapsto] & \Gamma_{\mQ}
\end{tikzcd}
\end{equation*}
in which the coset spaces for the three injections have canonical isomorphisms
\begin{equation}\label{eq: 5.2}
    \begin{tikzcd}
        \mathcal{G}_{\mQ}/\UGal \arrow{r}{\sim} & \cT_{\mQ}/\cT_F \arrow{r}{\sim} & \Gamma_{\mQ}/\Gamma_F \ar[r, "\sim"] & \mbox{Hom}_{\mQ}(F, \mC).
    \end{tikzcd}
\end{equation}
In particular, if $F/{\mQ}$ is a Galois extension, $\mathcal{G}_{\mQ}/\mathcal{G}_F$ is isomorphic to the Galois group $\Gamma_{{\mQ}/F}$.
\end{subconjecture}

This is Conjecture 4.2 of \cite{N2}. It is interesting in its own right, but it is also a step towards what we really want, a conjecture for the group over $\mQ$ whose complex points comprise our group $\UGal$. Another part of what we need is also available. It is the $\mQ$-structure defined by Langlands for the Taniyama group $\cT_F$ in the extension \eqref{eq: 4.6} that gives $\UGal$. What remains to account for is the complex group $\cD_F$ in \eqref{eq: 4.6}, and a $\mQ$-structure for the actual extension. These two conjectural steps were formulated briefly in \cite[Section 4]{N2}, with more details to be included in \cite{N14}. We shall be briefer still here.

The first step would be to define a quasisplit $F$-structure on $\cD_F$. This would begin with a quasisplit form $F$ on each of the factors $\cD_c$ in the product \eqref{eq: 4.7}. It would be obtained directly by transfer of the given quasisplit form of the group $G$ that is part of the associated index $c\sim (G, c)$. We would then compose this with an outer twist obtained by permutation of the factors $\{\cD_c\}$ by the action of the Galois group $\Gamma_F$ on the coordinates of the conjugacy classes $\{c_v : v\notin S\}$ in $\{\cD_c\}$. This is based on the understanding that for any $c$, the coordinates are defined over a finite extension $E = E_c$ of $F$. It also seems reasonable to suppose in the special case $F = \mQ$, that this gives the conjectured motivic $\mQ$-structure on the group $\mathcal{G} = \mathcal{G}_{\mQ}$.

Like the first step, the second can be described as in \cite[Section 4]{N2}. It is to expand the action of $\Gamma_F$ on $\cD_F$ to an action of $\Gamma_{\mQ}$ on the full group $\UGal$. This should include an inner twist that describes the dependence of $\UGal$, treated now as a group over $\mQ$ by restriction of scalars,
on the original fixed embedding $F\subset\mC$. We are assuming here as above that the required $\mQ$-structure on $\UGal$ is quasisplit, and defined by the first step with $F = \mQ$. For general fields $F\subset\mC$, we would need this to describe the required inner twists. I have not thought completely about this second step. Rough speculation asks whether it might be related to the inner twists attached by Langlands in Sections 6 and 7 of \cite{L4} to the group $G$ that defines a Shimura variety.

Our symbol $\UGal$ now has two interpretations. It represents either a projective limit of complex groups, as been the case until now, or a proalgebraic group $\cG = \cG_{\mQ}$ over $\mQ$. Of course the two interpretations are essentially the same. With its $\mQ$-structure, we restate Conjecture \ref{conj: 2(b)} for the associated groups of rational points as follows.

\begin{subconjecture}\label{conj: 2(c)}
    The groups of rational points $\mathcal{G}_F(\mathbb{Q})$ and $\cG(\mQ) = \mathcal{G}_{\mathbb{Q}}(\mathbb{Q})$ fit into a commutative diagram of prodiscrete homomorphisms
    \begin{equation*}
    \begin{tikzcd}
        \UGal(\mQ)\ar[r, mapsto]\ar[d, hook] & \cT_F(\mQ) \ar[r, mapsto]\ar[d, hook] & \Gamma_F\ar[d, hook]\\
        \mathcal{G}(\mQ)\ar[r, mapsto] & \cT(\mQ) \ar[r, mapsto] & \Gamma_Q,
\end{tikzcd}
\end{equation*}
in which the coset spaces for the three injections have canonical isomorphisms
\begin{equation}\label{eq: 5.3}
    \begin{tikzcd}
        \mathcal{G}(\mQ)/\UGal(\mQ) \arrow{r}{\sim} & \cT(\mQ)/\cT_F(\mQ) \arrow{r}{\sim} & \Gamma_Q/\Gamma_F \ar[r, "\sim"] & \mbox{Hom}_{\mQ}(F, \mC).
    \end{tikzcd}
\end{equation}
In particular, if $F/\mathbb{Q}$ is a Galois extension, $\mathcal{G}(\mathbb{Q})/\mathcal{G}_F(\mathbb{Q})$ is isomorphic to the Galois group $\Gamma_{\mQ/F}$.
\end{subconjecture}

This is Conjecture 4.3 of \cite{N2}. It relates to other articles that have studied $\mQ$-structures on automorphic representations. For the group $\GL(N)$ over $F$, Clozel introduced a property he called \ubf{algebraic} to describe what it should mean for an automorphic representation of this group to be motivic \cite{C}. Buzzard and Gee \cite{BG} expanded on this theme by studying analogues of Clozel's work for an arbitrary reductive group $G$. They defined four closely related properties that an automorphic representation $\pi$ of $G$ could have, which they call $C$-algebraic, $L$-algebraic, $C$-arithmetic and $L$-arithmetic, in terms of its local properties at certain valuations $\{v\}$, We discuss the questions raised in their paper in parts of Section 4 and 5 (and in a footnote for Section 6) of our paper \cite{N2}.

\section{Periods and more groups}\label{sec:periods and more groups}

The groups $L_F$ and $\UGal$, however they might appear now, are really just a beginning. They will have two kinds of natural extensions that are likely to dominate them in complexity (but also perhaps in richness!). The motivic groups will come with $\mQ$-algebras of complex numbers, called \ubf{periods}, on which the groups act in a way that extends the action of their quotient $\Gamma_{\mQ}$ on the corresponding subalgebras $\overline{\mQ}$. Natural extensions of the automorphic groups $L_F$ should also be expected, but there does not yet seem to be much information about them in the literature. We shall therefore confine our remarks to the natural extensions of $\UGal$. The two kinds of extensions would apply to \ubf{mixed} motives and \ubf{exponential} motives.

The motives we have discussed so far are often called \ubf{pure} motives. They are attached to nonsingular, projective algebraic varieties. But if we were to think of a \say{general} algebraic variety, as the set of complex valued solutions of several polynomial equations in a certain number of variables, it would not be projective. It would also be likely to be singular. However, there is an expanded theory of motives that applies to such general varieties. The objects so obtained are known as mixed motives. Possible references include \cite{Andre1} and \cite{H}.

According to general theory, largely conjectural, mixed motives also form a Tannakian category. With Betti cohomology over $\mQ$ again as a fibre functor, it corresponds to a proalgebraic group $\UGal^+$ over $\mQ$. In general, any algebraic group can be written uniquely as a semidirect product of a reductive group with a unipotent group. The proreductive part of the mixed motivic group would be the group $\UGal$ over $\mQ$ from the last section. The larger group would then be a semidirect product
\begin{equation*}
    \mathcal{G}^+_F = \mathcal{G}_F \ltimes \mathcal{U}_F,
\end{equation*}
for a unipotent proalgebraic group $\mathcal{U}_F$ over $\mathbb{Q}$. The goal would be to describe $\mathcal{U}_F$ in explicit terms.

The problem would seem in some ways to be within reach, even if there is not a conjectural solution at this point. Following the discussion of Section 9 of \cite{N2}, we note that as a unipotent group over $\mQ$, $\cU_F$ is completely determined by its Lie algebra  $\fu_F$. This in turn comes with its commutator quotient
\begin{equation*}
    \mathfrak{u}_{F, 0} = \mathfrak{u}_F/\mathfrak{u}_F^{der},
\end{equation*}
an infinite dimensional proalgebraic vector space. The action of $\mathcal{G}_F$ on $\cU_F$ transfers to a linear profinite dimensional representation of $\mathcal{G}_F$ on $\mathfrak{u}_F$ over $\mathbb{Q}$. On the other hand, it is expected that
\begin{equation}\label{eq: 6.1}
    \mbox{Ext}^i(M_1, M_2) = 0, \SP i\ge 2,
\end{equation}
for any pair $M_1$ and $M_2$ of mixed motives, according to one of the conjectures of Beilinson \cite[p.12]{N}. This means that $\mathcal{U}_F$ is a free unipotent group generated by the abelian quotient $\mathcal{U}_{F, 0} = \exp(\mathfrak{u}_{F, 0})$. In other words, if we understand $\mathfrak{u}_{F, 0}$ as a profinite direct sum of irreducible representations of $\mathcal{G}_F$ over $\mathbb{Q}$, we will in principle have a concrete description of $\mathcal{U}_{F}$ as a prounipotent group over $\mathbb{Q}$, armed with an action of $\mathcal{G}_F$. This would allow us in principle to extend our construction of $\mathcal{G}_F$ to a concrete description of the mixed motivic Galois group $\mathcal{G}^+_F$.

What then are the irreducible representations of $\UGal$, and for each of these, what is its multiplicity in $\fu_{F, 0}$? The constituents $\cD_c$ of $\UGal$ are simple, simply connected algebraic groups over $\mQ$. Their irreducible representations over $\mC$ would be classified by the Weyl theory of highest weights. This is a good start on the first question. One would then need to understand how these representations would all fit together to form irreducible representations of $\UGal$ over $\mQ$. Several things would be happening together at this point, one of which would concern the irreducible representations of the Taniyama group $\cT_F$. It is not clear how explicitly one would ultimately need to describe the irreducible representations of $\UGal$ over $\mQ$, or even how much is required in the description of their multiplicities. What does seem clear is that there is a good deal of structure to the problem, and that it would be very interesting to study.

The second question on the multiplicities is related to the work of F. Brown in Section $6$ of his paper \cite{B2c}\footnote{I thank Alexander Schmakov for pointing this out to me, and for informative conversations on the Beilinson conjectures.}, beginning with Proposition 6.1 there. Again, it is not clear how explicitly one would want to formulate the answer. In any case, using \eqref{eq: 6.1} one should then be able to describe $\cU_F$ as a prounipotent group over $\mQ$ with an action of $\UGal$, and hence $\UGal^+$ as a proalgebraic group over $\mQ$.

We far so far neglected one of the most fundamental properties of motives. It is the main reason for the adjective \say{Galois}. Both universal groups $L_F$ and $\UGal$ have the classical Galois groups $\Gamma_F$ as a quotient. However, the real analogy is in the totally disconnected group $\UGal(\mQ)$ of rational points, which behaves very much like a traditional Galois group. For there is a natural algebra $\cP_F(\mQ)$, the algebra of \ubf{$F$-periods}, with an action
\begin{equation*}
    \UGal(\mQ)\times \cP_F(\mQ) \longrightarrow \cP_F(\mQ)
\end{equation*}
of $\UGal(\mQ)$ that extends the action
\begin{equation*}
    \Gamma_F\times \overline{\mQ} \longrightarrow  \overline{\mQ}
\end{equation*}
of the field $\Gamma_F$ on $\overline{\mQ}$ of algebraic numbers. It contains $\overline{\mQ}$, but also many numbers that are expected to be transcendental. These typically take the form of definite integrals or convergent series that have long been part of classical analysis, but which could not be evaluated in elementary terms.

Grothendieck constructed periods in terms of his period isomorphism
\begin{equation}\label{eq: 6.2}
\begin{tikzcd}
    \varpi_X: H_{dR}^n(X)\otimes_F\mC \ar[r, "\sim"] & H_B^n(X)\otimes_{\mQ}(\mC),
\end{tikzcd}
\SP X\in \cP(F),
\end{equation}
where $H_{dR}^n(X)$ is algebraic de Rham cohomology and $H_B^n(X)$ is rational Betti cohomology, vector spaces over $F$ and $\mQ$ respectively. The mapping $\varpi_X$ is the complex isomorphism that represents the classical (complex) de Rham theorem for the compact, complex manifold $X(\mC)$. It is the matrix coefficients of $\varpi_X$ with respect to bases of $H_{dR}^n(X)$ and $H_B^n(X)$ (over $F$ and $\mQ$) that define periods. As $X$ varies, they generate the algebra $\cP_F(\mQ)$ of all periods. A reader unfamiliar with these notions could consult references \cite{H}, \cite[Section 7]{Andre1} or \cite[Section 7]{N2}.

The construction is conjectured to extend to arbitrary varieties/schemes $X$ over $F$. One would obtain an algebra $\PQ$ over $\mQ$ of mixed motivic periods, with an action
\begin{equation*}
    \UGal^+(\mQ) \times \PQ \longrightarrow \PQ
\end{equation*}
of the group of rational points in the mixed motivic Galois group $\UGal^+$. Examples of periods, mixed or pure, include the numbers $\pi$, $\log(q)$ for $q\in\mQ\backslash\{-1, 0, 1\}$, periods of elliptic curves over $F$, and special values of solutions of linear differential equations with regular singular points (and whose coefficients are rational functions with rational coefficients). They also include sums of many convergent infinite series. Notable in this last category are the values at integers $s\ge 2$ of the Riemann zeta function
\begin{equation*}
    \zeta(s) = \displaystyle\sum_{n = 1}^{\infty} (n^{-s}).
\end{equation*}
As is well known, its value $\zeta(2n)$ at an even integer is a rational multiple of $\pi^{2n}$, and hence a \ubf{pure} motivic period. On the other hand, its value $\zeta(2n + 1)$ at an odd integer has no such simple formula. However, it is still known to represent a \ubf{mixed} motivic period.

There is one number that is conspicuous for its absence from any such lists. It is the exponential constant $e$. This is not a (mixed) motivic period. It is rather a representative of a less familiar class of numbers, known (appropriately) as \ubf{exponential} periods. The theory of exponential periods is less developed, but it does come with a clear and comprehensive introduction \cite{FJ}, as well as very interesting work by various other authors. Like the theory of ordinary motivic periods, it is founded on general algebraic varieties over a number field $F\subset\mC$. What is different is its generalization of cohomology.

At the beginning of the introduction to \cite{FJ} the authors describe what they call \ubf{rapid decay homology}. It is attached to \ubf{varieties with potential}, by which is meant pairs $(X, f)$, where $X$ is a complex variety and $f:X \longrightarrow \mC$ is a regular function. It is defined as the limit
\begin{equation*}
    H_n^{rd}(X, f) = \lim_{r\longrightarrow\infty} H_n(X(\mC), f^{-1}(S_r); \mQ), \SP n\in \mathbb{N},
\end{equation*}
where the right hand side represents relative Betti homology with $\mQ$ coefficients, for the topological space $X(\mC)$ and the $f$-inverse images of closed half planes
\begin{equation*}
    S_r = \{z\in\mC \;:\; \mbox{Re}(z) \ge r\}.
\end{equation*}
With this, they are free to define rapid decay cohomology as a colimit
\begin{equation*}
    H_{rd}^n(X, f) = \mbox{Hom}_{\mQ}(H_n^{rd}(X, f), \mQ) = \underset{r\longrightarrow\infty}{\mbox{colim }} H^n(X(\mC), f^{-1}(S_r); \mQ).
\end{equation*}
(See \cite[Page 8]{FJ}). On the other hand, rapid decay cohomology has a purely algebraic counterpart, for $X$ smooth and defined over $F\subset\mC$. It is a generalization to $(X, f)$ of Grothendieck's algebraic de Rham cohomology \cite{G}, and gives a vector space $H^n_{dR}(X, f)$ over $F$.
The authors point out further that with a series of papers by various authors following an initial letter of Deligne \cite{N15}, one can construct a canonical isomorphism
\begin{equation}\label{eq: 6.3}
\begin{tikzcd}
    \varpi(X, f): H_{dR}^*(X, f)\otimes_F\mC \ar[r, "\sim"] & H_{rd}^*(X, f)\otimes_{\mQ}(\mC),
\end{tikzcd}
\end{equation}
between the two kinds of cohomology attached to pairs $(X, f)$. This reduces to the classical Grothendieck period formula \eqref{eq: 6.2} in the case that $f$ is a constant function. Moreover, as in its predecessor, the two sides of \eqref{eq: 6.3} extend, both separately and together as a formula, to arbitrary varieties over $F$. The matrix coefficients of the complex linear automorphism $\varpi_{X, f}$, with respect to bases of $H_{dR}^*(X, f)$ and $H_{rd}^*(X, f)$ (over $F$ and $\mQ$), are called \ubf{exponential} (mixed) \ubf{periods}. As $X$ varies, they should generate an algebra $\widetilde{\cP}_F^+(\mQ)$ over $\mQ$.

What we have taken here from the introduction to \cite{FJ} follows the inverse order of many basic presentations of motives and their periods. For we don't yet have exponential motives even though we have described their exponential periods. Motives are supposed to provide universal cohomology groups, from which all other formulations of cohomology are obtained as \ubf{realization} functors, together with comparison isomorphisms between them where relevant.

Because exponential periods come with two realizations $H_{dR}^*(X, f)$ and $H_{rd}^*(X, f)$ of de Rham cohomology, and an accompanying comparison isomorphism, it is suggested that we should expect a corresponding universal theory of exponential motives. This is indeed the case.

As is pointed out in \cite{FJ}, there is already a precedent. It is the theory of Nori for mixed motives, based on \ubf{quivers} and the \ubf{diagram category} among other things. Like much else, it relies on some basic conjectures. But with this proviso, it does provide a theory of mixed motives, \ubf{and} the mixed motivic Galois group $\UGal^+$.

Much of the volume \cite{FJ} is devoted to the extension of Nori's theory to exponential motives. The hardest part of this is the proof in Chapter 4 that the category of exponential motives is Tannakian. This in turn relies on a basic lemma of Beilinson on perverse sheaves. In the end, however, the category of exponential motives is constructed in \cite{FJ}. It again comes with a universal group, this time the exponential motivic Galois group $\widetilde{\mathcal{G}}_F^+$. We are thus presented with a variant of the question that has been the theme of this article, \say{Is there an explicit construction of the group $\widetilde{\mathcal{G}}_F^+$ that is founded on automorphic representations?} 

The exponential \ubf{pure} motivic Galois group should be given as an extension
\begin{equation}\label{eq: 6.4}
    1 \longrightarrow \widetilde{\mathcal{D}}_F\longrightarrow \widetilde{\mathcal{G}}_F \longrightarrow \widetilde{\mathcal{T}}_F \longrightarrow 1
\end{equation}
analogous to \eqref{eq: 4.6}. The group $\widetilde{\mathcal{T}}_F$ is known \cite{Ander}. It is given as an extension
\begin{equation*}
    1 \longrightarrow \widetilde{\mathcal{S}}_F\longrightarrow \widetilde{\mathcal{T}}_F \longrightarrow \Gamma_F \longrightarrow 1
\end{equation*}
analogous to \eqref{eq: 4.4}. Anderson referred to the objects attached to $\widetilde{\mathcal{T}}_F$ as \ubf{ulterior motives} in his paper. Fressan and Jossen showed in \cite{FJ} that they were indeed part of the general theory of exponential motives. The essential property of $\widetilde{\mathcal{G}}_F$ is that it maps surjectively onto $\mathcal{G}_F$, and hence fits into a commutative diagram
\begin{equation*}
  \begin{tikzcd}
      1 \ar[r] &\widetilde{\mathcal{D}}_F\ar[r]\ar[d]  &\widetilde{\mathcal{G}}_F \ar[r]\ar[d] &\widetilde{\mathcal{T}}_F \ar[r]\ar[d]  &1\\
      1 \ar[r] &\mathcal{D}_F\ar[r]  &\mathcal{G}_F\ar[r] &\mathcal{T}_F \ar[r]  &1.
  \end{tikzcd} 
\end{equation*}
The basic problem for us is to find an explicit construction $\widetilde{\mathcal{G}}_F$ like the one we described for $\UGal$. It would entail an automorphic version $\widetilde{L}_F$ of $L_F$, which among other things, would require an exponential generalization of the Shimura-Taniyama-Weil conjecture. In other words, we should be looking for a locally compact extension $\widetilde{L}_F$ of $L_F$, equipped with an extension $\widetilde{W}_F$ of the global Weil group, that embeds into a commutative diagram
\begin{equation*}
    \begin{tikzcd}
        1 \ar[r] & \widetilde{K}_F \ar[d] \ar[r] & \widetilde{L}_F \ar[d] \ar[r] & \widetilde{W}_F \ar[d]\ar[r] & 1\\
        1 \ar[r] &\widetilde{\mathcal{D}}_F\ar[r]  &\widetilde{\mathcal{G}}_F \ar[r] &\widetilde{\mathcal{T}}_F \ar[r]  &1\\
    \end{tikzcd}
\end{equation*}
analogous to \eqref{eq: 4.1}. What could it be? We would seem to have already used all of the relevant automorphic representations in our construction of $L_F$.

We conjecture that $\widetilde{L}_F$ will be defined by a construction similar to that of Section \ref{sec: automorphic galois groups}, but with \ubf{topological automorphic representations} in place of standard automorphic representations. These would be irreducible representations of topological covering groups
\begin{equation*}
    G\widetilde{(\mathbb{A}}) \longrightarrow G(\mathbb{A})
\end{equation*}
that occur in the decomposition of 
\begin{equation*}
    L^2(G(F) \backslash  G\widetilde{(\mathbb{A}})).
\end{equation*}
The Langlands program for topological automorphic representations is largely undeveloped, apart from exceptions for metaplectic representations, notably papers \cite{WL2}, \cite{WL1} by W-W.Li on metaplectic automorphic representations and the stabilization of associated trace formulas. In general, one would need a conjectural theory of Langlands dual groups, his principle of functoriality, and a resulting theory of $L$-functions. Present indications suggest that duals of groups $ G\widetilde{(\mathbb{A}})$ would not be coverings of groups $\widehat{G}$, but rather copies of (possibly different) classical dual groups. One could begin with what is known for metaplectic groups, in combination with the study of classical papers \cite{GGW}, \cite{MW} devoted to problems encountered for more general groups $ G\widetilde{(\mathbb{A}})$.

On the other hand, there has been considerable work on the actual extensions $ G\widetilde{(\mathbb{A}})$. These are not algebraic adele groups, but rather finite topological covering groups of topological groups $G(\mA)$, $G$ being a reductive group (say quasisplit) over $F$. The study of metaplectic groups, certain covers of order $2$, is classical. But there are many other topological extensions of general adelic groups $G(\mA)$ and their local factors $\widetilde{G}(F_v)$. Some of these are given by the theory of Brylinski-Deligne \cite{BD} that is founded on the algebraic $K$-theory group $K_2(F)$, but there are still others. If our conjecture holds, at least in principle, it will be very interesting to see which coverings are related through their representations to exponential motives. What seems clear is that topological automorphic representations greatly outnumber standard automorphic representations, just as we see from the definitions in \cite[Section 1.1.1]{FJ} that exponential motives will greatly outnumber classical motives.

The volume \cite{FJ} includes exponential \ubf{mixed} motives. There should then be an exponential mixed motivic Galois group $\widetilde{\mathcal{G}}_F^+$. Like the ordinary mixed motivic Galois group, it would be a proalgebraic group over $\mQ$, given as a semidirect product 
\begin{equation*}
    \widetilde{\mathcal{G}}_F^+ = \mathcal{G}_F \ltimes \widetilde{\mathcal{U}}_F
\end{equation*}
with its unipotent radical $\widetilde{\mathcal{U}}_F$ being a normal, proalgebraic unipotent group over $\mQ$.

If the conjecture above is correct, it opens the way for a major enhancement of the theory of automorphic forms. We know that there are many topological extensions of groups $G(\mA)$. There has been work on their representations, but as we noted above, a conjectural framework for how they extend the Langlands program is missing. Do they lead to new global $L$-functions? More likely is the possibility that classical Langlands $L$-functions come together in an expanded setting. We would then ask that the extended theory be parallel in some sense to that of exponential motives. What then would be the analogue of the Beilinson conjectures, for example? There are some indications that they bear upon values of automorphic $L$-functions at rational numbers $q>1$.

We are still missing something! If the extensions of the motivic Galois groups $\UGal$ make sense as described, we might expect them to be reflected in the automorphic Galois group $L_F$. If so, we would expect a mixed automorphic Galois group $L_F^+ = L \ltimes N_F$, an exponential automorphic Galois group $\widetilde{L}_F$ and an exponential mixed automorphic Galois group $\widetilde{L}_F^+ = L_f \ltimes \widehat{N}_F$. The best place to begin a search for such thing might be the (expanding) long article \cite{BSV}. This would be a long term commitment, which could be coordinated with the study of \ubf{automorphic periods}. Like motivic Galois groups, these objects might perhaps be paired naturally with our various enhancements of the automorphic Galois group.

To help keep track of the groups, we display the possible enhancements of $L_F$ and $\UGal$ and the maps among them in a commutative diagram as follows.
\begin{equation*}
    \begin{tikzcd}
       \Gamma_F \ar[r, phantom, "="] & \Gamma_F \ar[r, phantom, "="] & \Gamma_F \ar[r, phantom, "="] & \Gamma_F \ar[r, phantom, "="] & \Gamma_F \\
       W_F \ar[r, phantom, "="]\ar[u] & W_F \ar[u] & \widetilde{W}_F\ar[l]\ar[r, phantom, "="] \ar[u]& \widetilde{W}_F\ar[u]\ar[r] & W_F \ar[u]\\
        L_F \ar[u]\ar[d] & L_F^+\ar[u]\ar[d] \ar[l]& \widetilde{L}_F^+\ar[u]\ar[d]\ar[l] \ar{r} & \widetilde{L}_F\ar[u]\ar[d]\ar[r] & L_F\ar[u]\ar[d] \\
        \UGal \ar[d] & \UGal^+\ar[d]\ar[l] & \widetilde{\mathcal{G}}_F^+\ar[l]\ar[r]\ar[d] & \widetilde{\mathcal{G}}_F\ar[d]\ar[r] & \UGal \ar[d] \\
        \cT_F \ar[r, phantom, "="]\ar[d] & \cT_F\ar[d] & \widetilde{\cT}_F\ar[d]\ar[r]\ar[l] & \widetilde{\cT}_F\ar[d]\ar[r] & \cT_F\ar[d]\\
         \Gamma_F \ar[r, phantom, "="]& \Gamma_F \ar[r, phantom, "="] & \Gamma_F \ar[r, phantom, "="] & \Gamma_F  \ar[r, phantom, "="]& \Gamma_F
    \end{tikzcd}
\end{equation*}


\bibliographystyle{acm}
\bibliography{Bibliography}

\begin{thebibliography}{10}

\bibitem{Ander}
{\sc Anderson, G.~W.}
\newblock Cyclotomy and an extension of the {T}aniyama group.
\newblock {\em Compositio Math. 57}, 2 (1986), 153--217.

\bibitem{Andre1}
{\sc Andr{\'e}, Y.}
\newblock Une introduction aux motifs (motifs purs, motifs mixtes,
  p{\'e}riodes).
\newblock In {\em Panoramas et Synthèses\/} (2004).

\bibitem{Ar2}
{\sc Arthur, J.}
\newblock A note on the automorphic {L}anglands group.
\newblock vol.~45. 2002, pp.~466--482.
\newblock Dedicated to Robert V.\ Moody.

\bibitem{Ar3}
{\sc Arthur, J.}
\newblock {\em The endoscopic classification of representations}, vol.~61 of
  {\em American Mathematical Society Colloquium Publications}.
\newblock American Mathematical Society, Providence, RI, 2013.
\newblock Orthogonal and symplectic groups.

\bibitem{N14}
{\sc Arthur, J.~G.}
\newblock Structure on universal groups.
\newblock In preparation.

\bibitem{N1}
{\sc Arthur, J.~G.}
\newblock The work of {R}obert {L}anglands.
\newblock In {\em The {A}bel {P}rize 2018--2022}. Springer, Cham, 2024,
  pp.~31--230.

\bibitem{N2}
{\sc Arthur, J.~G.}
\newblock On {U}niversal {G}alois {G}roups, 2025.
\newblock To appear.

\bibitem{AGIKMS}
{\sc Atobe, H., Gan, W.~T., Ichino, A., Kaletha, T., Mínguez, A., and Shin,
  S.~W.}
\newblock Local intertwining relations and co-tempered $a$-packets of classical
  groups, 2024.

\bibitem{BSV}
{\sc Ben-Zvi, D., Sakellaridis, Y., and Venkatesh, A.}
\newblock Relative langlands duality, 2024.

\bibitem{N4}
{\sc Borel, A., and Harish-Chandra}.
\newblock Arithmetic subgroups of algebraic groups.
\newblock {\em Bull. Amer. Math. Soc. 67\/} (1961), 579--583.

\bibitem{BJ}
{\sc Borel, A., and Jaquet, H.}
\newblock Automorphic forms and automorphic representations.
\newblock In {\em Automorphic Forms, Representations and L-functions}, vol.~33
  part 1 of {\em Symposium in Pure Mathematics}. American Mathematical Society,
  1979, pp.~189--202.

\bibitem{B2c}
{\sc Brown, F.}
\newblock Notes on motivic periods.
\newblock {\em Commun. Number Theory Phys. 11}, 3 (2017), 557--655.

\bibitem{BD}
{\sc Brylinski, J.-L., and Deligne, P.}
\newblock Central extensions of reductive groups by {$\bold K_2$}.
\newblock {\em Publ. Math. Inst. Hautes \'Etudes Sci.}, 94 (2001), 5--85.

\bibitem{BG}
{\sc Buzzard, K., and Gee, T.}
\newblock The conjectural connections between automorphic representations and
  {G}alois representations.
\newblock In {\em Automorphic forms and {G}alois representations. {V}ol. 1},
  vol.~414 of {\em London Math. Soc. Lecture Note Ser.} Cambridge Univ. Press,
  Cambridge, 2014, pp.~135--187.

\bibitem{C}
{\sc Clozel, L.}
\newblock Motifs et formes automorphes: applications du principe de
  fonctorialit\'e.
\newblock In {\em Automorphic forms, {S}himura varieties, and {$L$}-functions,
  {V}ol.\ {I} ({A}nn {A}rbor, {MI}, 1988)}, vol.~10 of {\em Perspect. Math.}
  Academic Press, Boston, MA, 1990, pp.~77--159.

\bibitem{CR}
{\sc Clozel, L., and Rajan, C.}
\newblock Solvable base change.
\newblock {\em Journal für die reine und angewandte Mathematik\/} (06 2020).

\bibitem{N15}
{\sc Deligne, P.}
\newblock Letter to {M}algrange, 1976.

\bibitem{FJ}
{\sc Fresan, J., and Jossen, P.}
\newblock Exponential motives.
\newblock preprint.

\bibitem{GGW}
{\sc Gan, W.~T., Gao, F., and Weissman, M.~H.}
\newblock L-groups and the {L}anglands program for covering groups: a
  historical introduction, 2017.

\bibitem{N10}
{\sc Godement, R., and Jacquet, H.}
\newblock {\em Zeta functions of simple algebras}, vol.~Vol. 260 of {\em
  Lecture Notes in Mathematics}.
\newblock Springer-Verlag, Berlin-New York, 1972.

\bibitem{G}
{\sc Grothendieck, A.}
\newblock Standard conjectures on algebraic cycles.
\newblock In {\em Algebraic {G}eometry ({I}nternat. {C}olloq., {T}ata {I}nst.
  {F}und. {R}es., {B}ombay, 1968)}, vol.~4 of {\em Tata Inst. Fundam. Res.
  Stud. Math.} Tata Inst. Fund. Res., Bombay, 1969, pp.~193--199.

\bibitem{N3}
{\sc Harish-Chandra}.
\newblock Automorphic forms on a semisimple {L}ie group.
\newblock {\em Proc. Nat. Acad. Sci. U.S.A. 45\/} (1959), 570--573.

\bibitem{N7}
{\sc Hecke, E.}
\newblock Eine neue {A}rt von {Z}etafunktionen und ihre {B}eziehungen zur
  {V}erteilung der {P}rimzahlen.
\newblock {\em Math. Z. 6}, 1-2 (1920), 11--51.

\bibitem{N8}
{\sc Hecke, E.}
\newblock \"uber {M}odulfunktionen und die {D}irichletschen {R}eihen mit
  {E}ulerscher {P}roduktentwicklung. {I}.
\newblock {\em Math. Ann. 114}, 1 (1937), 1--28.

\bibitem{H}
{\sc Huber, A., and M\"uller-Stach, S.}
\newblock {\em Periods and {N}ori motives}, vol.~65 of {\em Ergebnisse der
  Mathematik und ihrer Grenzgebiete. 3. Folge. A Series of Modern Surveys in
  Mathematics [Results in Mathematics and Related Areas. 3rd Series. A Series
  of Modern Surveys in Mathematics]}.
\newblock Springer, Cham, 2017.
\newblock With contributions by Benjamin Friedrich and Jonas von Wangenheim.

\bibitem{N9}
{\sc Jacquet, H., and Langlands, R.~P.}
\newblock {\em Automorphic forms on {${\rm GL}(2)$}}, vol.~Vol. 114 of {\em
  Lecture Notes in Mathematics}.
\newblock Springer-Verlag, Berlin-New York, 1970.

\bibitem{KO}
{\sc Kottwitz, R.~E.}
\newblock Stable trace formula: cuspidal tempered terms.
\newblock {\em Duke Math. J. 51}, 3 (1984), 611--650.

\bibitem{N6}
{\sc langlands, R.}
\newblock Letter to {A}. {W}eil, 1967.

\bibitem{L2}
{\sc Langlands, R.~P.}
\newblock Problems in the theory of automorphic forms.
\newblock In {\em Lectures in {M}odern {A}nalysis and {A}pplications, {III}},
  vol.~Vol. 170 of {\em Lecture Notes in Math.} Springer, Berlin-New York,
  1970, pp.~18--61.

\bibitem{N5}
{\sc Langlands, R.~P.}
\newblock {\em On the functional equations satisfied by {E}isenstein series},
  vol.~Vol. 544 of {\em Lecture Notes in Mathematics}.
\newblock Springer-Verlag, Berlin-New York, 1976.

\bibitem{L4}
{\sc Langlands, R.~P.}
\newblock Automorphic representations, {S}himura varieties, and motives. {E}in
  {M}\"archen.
\newblock In {\em Automorphic forms, representations and {$L$}-functions
  ({P}roc. {S}ympos. {P}ure {M}ath., {O}regon {S}tate {U}niv., {C}orvallis,
  {O}re., 1977), {P}art 2}, vol.~XXXIII of {\em Proc. Sympos. Pure Math.} Amer.
  Math. Soc., Providence, RI, 1979, pp.~205--246.

\bibitem{L5}
{\sc Langlands, R.~P.}
\newblock On the notion of an automorphic representation: A supplement to the
  proceeding paper.
\newblock In {\em Automorphic forms, representations and {$L$}-functions
  ({P}roc. {S}ympos. {P}ure {M}ath., {O}regon {S}tate {U}niv., {C}orvallis,
  {O}re., 1977), {P}art 2}, vol.~XXXIII of {\em Proc. Sympos. Pure Math.} Amer.
  Math. Soc., Providence, RI, 1979, pp.~203--208.

\bibitem{N12}
{\sc Langlands, R.~P.}
\newblock {\em Les d\'ebuts d'une formule des traces stable}, vol.~13 of {\em
  Publications Math\'ematiques de l'Universit\'e{} Paris VII [Mathematical
  Publications of the University of Paris VII]}.
\newblock Universit\'e{} de Paris VII, U.E.R. de Math\'ematiques, Paris, 1983.

\bibitem{N13}
{\sc Langlands, R.~P.}
\newblock On the classification of irreducible representations of real
  algebraic groups.
\newblock In {\em Representation theory and harmonic analysis on semisimple
  {L}ie groups}, vol.~31 of {\em Math. Surveys Monogr.} Amer. Math. Soc.,
  Providence, RI, 1989, pp.~101--170.

\bibitem{WL2}
{\sc Li, W.-W.}
\newblock Arthur packets for metaplectic groups, 2024.

\bibitem{WL1}
{\sc Li, W.-W.}
\newblock Stabilization of the trace formula for metaplectic groups, 2024.

\bibitem{Mo}
{\sc Mok, C.~P.}
\newblock Endoscopic classification of representations of quasi-split unitary
  groups.
\newblock {\em Mem. Amer. Math. Soc. 235}, 1108 (2015), vi+248.

\bibitem{N}
{\sc Nekov\'a\v~r, J.}
\newblock Be\u ilinson's conjectures.
\newblock In {\em Motives ({S}eattle, {WA}, 1991)}, vol.~55, Part 1 of {\em
  Proc. Sympos. Pure Math.} Amer. Math. Soc., Providence, RI, 1994,
  pp.~537--570.

\bibitem{R}
{\sc Saavedra~Rivano, N.}
\newblock {\em Cat\'egories {T}annakiennes}, vol.~Vol. 265 of {\em Lecture
  Notes in Mathematics}.
\newblock Springer-Verlag, Berlin-New York, 1972.

\bibitem{Se1}
{\sc Serre, J.-P.}
\newblock {\em Abelian {$l$}-adic representations and elliptic curves},
  second~ed.
\newblock Advanced Book Classics. Addison-Wesley Publishing Company, Advanced
  Book Program, Redwood City, CA, 1989.
\newblock With the collaboration of Willem Kuyk and John Labute.

\bibitem{Se2}
{\sc Serre, J.-P.}
\newblock Propri\'et\'es conjecturales des groupes de {G}alois motiviques et
  des repr\'esentations {$l$}-adiques.
\newblock In {\em Motives ({S}eattle, {WA}, 1991)}, vol.~55, Part 1 of {\em
  Proc. Sympos. Pure Math.} Amer. Math. Soc., Providence, RI, 1994,
  pp.~377--400.

\bibitem{T}
{\sc Tate, J.}
\newblock Number theoretic background.
\newblock In {\em Automorphic forms, representations and {$L$}-functions
  ({P}roc. {S}ympos. {P}ure {M}ath., {O}regon {S}tate {U}niv., {C}orvallis,
  {O}re., 1977), {P}art 2}, vol.~XXXIII of {\em Proc. Sympos. Pure Math.} Amer.
  Math. Soc., Providence, RI, 1979, pp.~3--26.

\bibitem{N11}
{\sc Weil, A.}
\newblock Sur la th\'eorie du corps de classes.
\newblock In {\em S\'eminaire {B}ourbaki, {V}ol.\ 2}. Soc. Math. France, Paris,
  1995, pp.~Exp. No. 83, 313--315.

\bibitem{MW}
{\sc Weissman, M.~H.}
\newblock A comparison of {L}-groups for covers of split reductive groups.
\newblock No.~398. 2018, pp.~277--286.
\newblock L-groups and the Langlands program for covering groups.

\end{thebibliography}

\end{document}